\documentclass[12pt]{article}
\usepackage{amssymb}

\newtheorem{thm}     {Theorem}[section]
\newtheorem{prop}    [thm]{Proposition}
\newtheorem{definition}  [thm]{Definition}
\newtheorem{cor}     [thm]{Corollary}
\newtheorem{lemma}   [thm]{Lemma}

\newcommand{\proof} {\noindent{\bf Proof. }}
\newcommand{\A}{{\mathcal A}}
\newcommand{\B}{\mathbb B}
\newcommand{\C}{\mathbb C}
\newcommand{\D}{\mathbb D}
\newcommand{\R}{\mathbb R}
\newcommand{\T}{\mathbb T}
\newcommand{\Z}{\mathbb Z}
\newcommand{\st}{{\rm st}}
\newcommand{\rh}{{\rm rh}}
\def\Re{{\rm Re\,}}

\def\bar{\overline}

\begin{document}

\title{Hartogs figure and symplectic non-squeezing}
\author{Alexandre Sukhov{*} and Alexander Tumanov{**}}
\date{}
\maketitle

\rightline{\it To John D'Angelo}
\bigskip
\bigskip
\bigskip

{\small
* Universit\'e des Sciences et Technologies de Lille, Laboratoire
Paul Painlev\'e,
U.F.R. de
Math\'e-matique, 59655 Villeneuve d'Ascq, Cedex, France,
 sukhov@math.univ-lille1.fr

** University of Illinois, Department of Mathematics
1409 West Green Street, Urbana, IL 61801, USA, tumanov@math.uiuc.edu
}
\bigskip

Abstract.
We solve a problem on filling by Levi-flat hypersurfaces
for a class of totally real 2-tori in a real 4-manifold with an
almost complex structure tamed by an exact symplectic form.
As an application we obtain a simple proof of Gromov's
non-squeezing theorem in dimension 4 and new results on
rigidity of symplectic structures.
\bigskip


\bigskip

\section{Introduction}

Since Gromov's work \cite{Gr} it is known that $J$-complex curves
can be used in order to describe obstructions for symplectic
embeddings. Following this theme, in this paper we apply classical
complex analysis to symplectic rigidity. We obtain new results on
non-existence of certain symplectic embeddings, in particular,
we give a simple proof of Gromov's non-squeezing theorem
in complex dimension 2.
Our approach is based on a general result on Levi-flat fillings
of totally real tori in an almost complex manifold with an exact
symplectic form. This result is new even for manifolds with
integrable almost complex structure.

\begin{definition}
\label{radius}
Let $G$ be a domain in $\C^2$ containing the origin.
Denote by ${\cal O}^1_0(G)$ the set of closed complex purely
one-dimensional analytic subsets in $G$ passing through the origin.
Denote by $E(X)$ the Euclidean area of $X \in {\cal O}^1_0(G)$.
The holomorphic radius $\rh(G)$ of $G$ is defined as
$$
\rh(G)= \inf \{ \lambda > 0: \exists X \in
{\cal O}^1_0(G), E(X) = \pi \lambda^2 \}.
$$
If the set in the right-hand part is empty, then we set
$\rh(G) = +\infty$.
\end{definition}

\noindent
{\bf Example.}
Let $\B$ be the Euclidean ball of $\C^2$
and $r> 0$. Then $\rh(r\B) = r$.
Indeed, the area $E(X)$ of $X \in {\cal O}^1_0(r\B)$ is bounded
from below by the area $\pi r^2$ of a section of the ball by
a complex line through the origin
(Lelong, 1950; see \cite{Ch}).
\bigskip

Let $(z_1,z_2)$, $z_j = x_j + iy_j$, be complex coordinates
in $\C^2$. Let
$$
\omega_\st = \frac{i}{2}\sum_{j=1}^2 dz_j \wedge d\bar{z}_j
$$
be the standard symplectic form on $\C^2$.
A diffeomorphism $\phi: G_1 \to G_2$ between two domains
$G_j\subset\C^2$ is called a {\it symplectomorphism} if
$\phi^*\omega_\st = \omega_\st$.
Let $\D$ denote the unit disc in $\C$.
Our main result concerning
symplectic rigidity is the following

\begin{thm}
\label{rigidity}
Let $G_1$ be a domain in $\C^2$ containing the origin
and let $G_2$ be a domain  in $R\D \times \C$ for some $R > 0$.
Assume that there exists a symplectomorphism
$\phi: G_1 \to G_2$. Then $\rh(G_1) \leq R$.
\end{thm}

In view of the example above we obtain Gromov's \cite{Gr}
non-squeezing theorem  in $\C^2$.
\begin{cor}
\label{Gromov}
Suppose that there exists a symplectomorphism between the
ball $r\B$ and a domain  contained in $R\D \times \C$.
Then $r \leq R$.
\end{cor}

As an alternative to the usual {\it complex bidisc}
$$
\D^2 = \{ (z_1,z_2):   |z_j|< 1, j=1,2 \}
$$
we introduce the {\it real bidisc} of the form
$$
\D^2_\R = \{ (z_1,z_2): x_1^2 + x_2^2 < 1, y_1^2 + y_2^2 < 1 \}.
$$
The two bidiscs have the same volume. Are they symplectomorphic?
We learned about this question from Sergey Ivashkovich \cite{Iv}.
We will prove that the answer is negative.
Note that if a symplectomorphism $\D^2\to\D^2_\R$ is smooth
up to the boundary, then it maps the torus
$\T^2=\{ (z_1,z_2):   |z_j|= 1, j=1,2 \}$ to the torus
$\T^2_\R=\{ (z_1,z_2): x_1^2 + x_2^2=1, y_1^2 + y_2^2=1 \}$.
However, it is not possible because $\T^2$ is Lagrangian,
that is, $\omega_\st|_{\T^2}=0$, while $\T^2_\R$ is not.
This may lead one to a thought that the question is about
exotic non-smooth maps. We show it is not the case.
In fact $\D^2_\R$ does not admit a symplectic embedding
into a slightly larger complex bidisc.
Furthermore, we obtain the following non-squeezing result.

\begin{cor}
\label{bidisc1}
There exists $R > 1$ such that there is no symplectomorphism
between $\D^2_\R$ and a subdomain of $R\D \times \C$.
\end{cor}

We will show in the last section that $\rh(\D^2_\R) > 1$, then
Corollary \ref{bidisc1} will immediately follow from
Theorem \ref{rigidity}. In particular, we obtain

\begin{cor}
\label{RealSqueez}
There is no symplectomorphism between the real bidisc $\D^2_\R$
and the complex bidisc $\D^2$.
\end{cor}

The proof of Theorem \ref{rigidity} relies on filling by complex
discs an analog of the Hartogs figure for an almost complex manifold.
Let $(M,J,\omega)$ be a $C^\infty$-smooth real 4-dimensional
manifold with a symplectic form $\omega$ and an almost complex
structure $J$. We suppose that $J$ is {\it tamed} by $\omega$
(see \cite{Gr}), i. e., $\omega(V,JV) > 0$ for every non-zero
tangent vector $V$; we call such  $M$
a {\it tame} almost complex manifold.
Consider  a relatively compact subdomain  $\Omega$ in $M$ with
smooth {\it strictly pseudoconvex boundary}. This means that for
every point $p \in b\Omega$ there exists an open neighborhood $U$
and a smooth strictly $J$-plurisubharmonic function
$\rho:U \to \R$ with non-vanishing gradient such that
$\Omega \cap U = \{ q \in U: \rho(q) < 0 \}$.
We do {\it not} require the existence of a global defining
strictly plurisubharmonic function on $\bar\Omega$.

We now use the notation $Z = (z,w)$ for complex coordinates
in $\C^2$.

\begin{definition}
\label{Hartogs}A $C^\infty$-smooth embedding
$H: \overline\D \times \overline \D\to \bar\Omega$
is called a Hartogs embedding if the following conditions hold:
\begin{itemize}
\item[(i)]
the map $f^0 : \D\ni z \mapsto H(z,0)$
is $J$-complex and $f^0(\bar\D) \subset \Omega$;
\item[(ii)]
for every $z \in \overline\D$, the map
$h_z : \D\ni w \mapsto H(z,w)\in \Omega$ is
$J$-complex; moreover, there exists $\delta > 0$ such that
for every $z$ with $1-\delta \le |z| \le 1$, we have
$h_z(b\D) \subset b\Omega$;
\end{itemize}
\end{definition}

Denote $\Lambda^t= H(b\D \times tb\D)$, $0 \leq t \leq 1$.
Then $\Lambda^t$ is a totally real torus in $M$.
We will also denote by $\Pi$ the Levi-flat hypersurface
$\Pi = H(b\D \times \D)$.
Thus the family of tori $\Lambda^t$ and the hypersurface
$\Pi$ are canonically associated with a Hartogs embedding.

Our main technical tool is the following

\begin{thm}
\label{tori1}
Let $\Omega$ be a relatively compact domain with smooth strictly
pseudoconvex boundary in a tame almost complex manifold
$(M,J,\omega)$ of complex dimension 2 and let
$H:\bar\D \times \bar\D \to \bar\Omega$ be a Hartogs embedding.
Assume that the symplectic form $\omega$ is exact in a neighborhood
of the closure $\bar \Omega$. Then for every $0 < t \leq 1$
there exists a unique one-parameter family of  embedded
$J$-complex discs $f:\D \to \Omega$ of class
$C^\infty(\overline\D)$ such that $f(b\D) \subset \Lambda^t$.
They fill a smooth Levi-flat hypersurface $\Gamma^t \subset \Omega$
with boundary $\Lambda^t$. The family $(\Gamma^t)$ foliates
a subdomain in $\Omega$ whose boundary consists of the Levi-flat
hypersurfaces $\Gamma^1$ and $\Pi$ and the disc $f^0(\bar\D)$.
\end{thm}

For simplicity we assume that $(M,J,\omega)$ is $C^\infty$,
however the proof needs a finite smoothness.
We construct the desired discs by a continuous deformation
starting from the initial disc $f^0$.
In particular, they are
homotopic to $f^0$ in the space of $J$-complex discs in
$\bar\Omega$ attached to $\Pi$. A similar approach was used
by Bedford and Gaveau \cite{BeGa}, Forstneri\v c \cite{Fo},
Gromov \cite{Gr}, and others in various situations.
The statement of Theorem \ref{tori1} can be slightly improved
by introducing the map $H$ with the given properties
only on $b\D \times \D$. We keep the stated version
for simplicity and convenience of presentation.

In his celebrated paper, Gromov \cite{Gr} proved that for every
compact {\it Lagrangian} submanifold $\Lambda$ in $\C^n$,
there exists a non-constant complex disc with boundary in $\Lambda$.
In comparison, our Theorem \ref{tori1} applies to non-Lagrangian
tori and it gives information about the set swept out
by the discs.
In the case $M=\C^2$ with the standard complex structure,
there are related results due to Duval-Gayet \cite{DuGa}
and Forstneri\v c \cite{Fo}.
We stress that Theorem \ref{tori1} is new even in the case
$M$ is a complex manifold, i. e., the structure $J$ is integrable.

Recall that the classical Hartogs figure $U$ is a neighborhood
of $(\D \times \{ 0 \}) \cup (b\D \times \D)$ in $\C^2$.
One can choose $U$ as a union of complex
discs $\{ z \} \times r(z)\D$, $z \in \bar\D$, where
$0 < r(z) \leq 1$ is smooth in $z$; then the embedding
$H: \bar\D^2 \to U$ defined by $H(z,w) = (z,r(z)w)$ satisfies
Definition \ref{Hartogs} and smoothly parametrizes the Hartogs
figure by the bidisc. Therefore, one can view Theorem \ref{tori1}
as a result on filling a Hartogs figure by complex discs.
In the case $J$ is integrable it can be used in the study
of holomorphic extension problems and  polynomial,
holomorphic, and plurisubharmonic hulls. We also point out that
the Hartogs embedding in Definition \ref{Hartogs} is not
necessarily biholomorphic, which brings additional flexibility
to the method. In this paper we focus on symplectic applications
of the theorem.

This paper was written for a special volume in honor of
our dear colleague Professor John D'Angelo on the occasion
of his 60th birthday. The authors wish John good health,
happiness, and new research accomplishments for years to come.

\section{Almost complex manifolds}

Let $(M,J)$ be an almost complex manifold.
We denote  by $J_\st$ the standard complex structure of $\C^n$;
the value of $n$ will be clear from the context.
A $C^1$-map $f:\D \to M$
is called a {\it $J$-complex} (or $J$-holomorphic) disc
if $df \circ J_\st = J \circ df$.

In local coordinates $Z=(z,w)\in\C^2$, an almost complex
structure $J$ can be represented by a complex $2\times 2$
matrix function $A$, so that a map $Z:\D \to\C^2$ is
$J$-complex if and only if it satisfies the following
partial differential equation
\begin{eqnarray}
\label{CR}
Z_{\overline\zeta} - A(Z)\overline{Z_\zeta}=0.
\end{eqnarray}
The matrix $A(Z)$ is defined by
\begin{eqnarray}
\label{matrixA} A(Z)V = (J_\st + J(Z))^{-1}
(J_\st -J(Z))\overline V.
\end{eqnarray}
Indeed, one can see that the right-hand side of (\ref{matrixA})
is $\C$-linear in $V\in\C^n$ with respect to the standard
structure $J_\st$, hence $A(Z)$ is well defined
(see, e. g., \cite{SuTu1}).
We call $A$ {\it the complex matrix} of $J$.
The ellipticity of (\ref{CR}) is equivalent to
$\det(I-A\bar A)\ne0$.
In a fixed coordinate chart, the  correspondence
between almost complex structures $J$ with $\det(J_\st+J)\ne0$
and complex matrices with $\det(I-A\bar A)\ne0$ is one-to-one
\cite{SuTu1}.

Often we identify $J$-complex discs  $f$ and
their images calling them just the discs. By the boundary
of such a disc we mean the restriction $f\vert_{b\D}$,
which we also identify with its image.

Let $\rho$ be a  function of class $C^2$ on $M$, let $p \in M$
and $V\in {T}_pM$. {\it The Levi form} of $\rho$ at $p$ evaluated
on $V$ is defined by the equality
$L^J(\rho)(p)(V):=-d(J^* d\rho)(V,JV)(p)$.
A  real function $\rho$ of class $C^2$ on $M$ is called
{\it  $J$-plurisubharmonic} (resp. strictly $J$-plurisubharmonic)
if  $ L^J(\rho)(p)(V) \geq 0$ (resp. $> 0$) for every $p \in M$,
$V \in {T}_pM \backslash \{0\}$.

A smooth real hypersurface $\Gamma$ in an almost complex
manifold $(M,J)$ is called {\it Levi-flat} if in a neighborhood
$U$ of  every point $p \in \Gamma$ there exists a defining
function with non-zero gradient whose Levi form vanishes
for every tangent vector
$V \in {T}_q\Gamma \cap J {T}_q\Gamma$
and every point $q \in U \cap \Gamma$.
If the complex dimension of $M$ is equal to $2$, then
by the Frobenius theorem, a hypersurface $\Gamma$ is Levi-flat
if and only if $\Gamma$ is locally foliated by a
real one-parameter family of $J$-complex discs.

\section{Deformation}

Returning to Theorem \ref{tori1}, we fix $\delta >0$
that figures in Definition \ref{Hartogs}.
Consider the annulus
$$
\A_\delta=\{z\in\C:1-\delta\le|z|\le 1\}.
$$
Introduce also the discs and the circles
$$
G^t=\{w\in\D:  \vert w \vert<t \},\qquad
\gamma^t=b G^t = \{ w \in \D: \vert w \vert = t \}.
$$
We recall the notations $\Lambda=H(b\D \times b\D)$,
$\Pi = H(b\D \times \bar\D)$
and
$$
\Lambda^t=\bigcup_{z\in b\D}
H(\{z\}\times\gamma^t), \quad
0\le t\le 1.
$$
Then for $0< t\le 1$, $\Lambda^t$ is a totally real torus,
$\Lambda^1=\Lambda$, and $\Lambda^0=f^0(b\D)$ is a circle.

We will consider $J$-complex discs with boundaries in
$\Lambda^t$. By reflection principle \cite{IvSu} such discs
are smooth up to the boundary.

Let $t_0>0$.
Let $I(t_0)$ denote one of the intervals: $[0,t_0]$ or $[0,t_0)$.
Let
$$
\{f^{t,\tau}:\bar\D\to M: t\in I(t_0), \tau\in\R/2\pi\Z\}
$$
be a continuous family of embedded
$J$-complex discs, smooth in all the variables for $t>0$.
\begin{definition}
\label{def-def}
We call the family $(f^{t,\tau})$
{\it an admissible deformation} (of the initial disc $f^0$)
on $I(t_0)$ if it has the following properties.
\begin{itemize}
\item[(i)]
$f^{0,\tau}=f^0$.

\item[(ii)]
$f^{t,\tau}(b\D)\subset\Lambda^t$;
$f^{t,\tau_1}(b\D)\cap f^{t,\tau_2}(b\D)=\emptyset$
if $\tau_1\ne\tau_2$;
$\bigcup_{\tau\in\R/2\pi\Z} f^{t,\tau}(b\D)=\Lambda^t$;
$\Gamma^t=\bigcup_{\tau\in\R/2\pi\Z}f^{t,\tau}(\bar\D)$
is a smooth hypersurface with boundary $\Lambda^t$.

\item[(iii)]
The set $H^{-1}(f^{t,\tau}(\bar\D))\cap (\A_\delta\times\bar\D)$
is the graph of a non-vanishing smooth function
$w^{t,\tau}:{\cal A}_\delta\to\bar\D \backslash \{ 0 \}$
so that the map $w^{t,\tau}|_{b\D}:b\D\to\gamma^t$
has zero winding number.
Furthermore,
$f^{t,\tau}(\bar\D)\cap f^0(\bar\D)=\emptyset$
for $t>0$.

\item[(iv)]
(normalization condition)
For a fixed $\zeta_0$ in the interior of $\A_\delta$,
say $\zeta_0=1-\delta/2$, we have $w^{t,\tau}(1)=te^{i\tau}$,
$f^{t,\tau}(1)=H(1,w^{t,\tau}(1))$, and
$f^{t,\tau}(\zeta_0)=H(\zeta_0,w^{t,\tau}(\zeta_0))$.

\item[(v)]
Every $J$-complex disc $f$ such that $f(b\D)\subset\Lambda^t$
and close to $f^{t,\tau}$ in $C^{1,\alpha}(\D)$,
coincides with $f^{t,\tau'}$ for some $\tau'\in\R/2\pi\Z$
close to $\tau$ after a reparametrization close to
the identity; here $0<\alpha<1$ is fixed, say $\alpha=1/2$.
(Note that $f\in C^\infty(\bar\D)$ by reflection principle
\cite{IvSu}.)
\end{itemize}
\end{definition}

Since $b\Omega$ is strictly pseudoconvex, for every $t$ and $\tau$
the discs $f^{t,\tau}$ are contained in $\Omega$.

Consider the pull-back $H^*(J)$ of $J$ to $\D^2$.
It follows from  Definition \ref{Hartogs} (see, e. g.,
\cite{SuTu2}) that the complex matrix $A$ of $H^*(J)$
over ${\cal A}_\delta \times \D$
has the following special form:
\begin{eqnarray}
\label{structure}
A= \left(
\begin{array}{cll}
a & & 0\\
b & & 0
\end{array}
\right)
\end{eqnarray}
with $\vert a \vert < 1$.
Then (see, e. g., \cite{SuTu2})
the functions $w^{t,\tau}$ satisfy the equation
\begin{eqnarray}
\label{CRgraph}
w_{\bar z} + a w_z = b
\end{eqnarray}
Conversely, the graph of every solution of (\ref{CRgraph})
becomes a $J$-complex curve after a suitable reparametrization
$z = z(\zeta)$.

We prove Theorem \ref{tori1} by showing that there
exists an admissible deformation on $[0,1]$.
Sometimes we will write $f^{t}$ instead of $f^{t,\tau}$
if the value of $\tau$ is unimportant.

In \cite{SuTu2} we obtain a result similar to Theorem \ref{tori1}
for $M=\C^2$ equipped with an almost complex structure
whose matrix has a form even more general than (\ref{structure}).
We can use that result in a neighborhood of the disc $f^0$.
Then we obtain the following
\begin{prop}
\label{deformation0}
For small $t_0>0$ there exists a unique admissible
deformation on $I(t_0)$.
\end{prop}

\begin{prop}
\label{Levi-flat1}
\begin{itemize}
\item[(i)]
If an admissible deformation $f^{t,\tau}$ on $I(t_0)$ exists,
then it is unique.
\item[(ii)]
$f^{t_1,\tau_1}(\bar\D)\cap f^{t_2,\tau_2}(\bar\D)=\emptyset$
unless $t_1=t_2$ and $\tau_1=\tau_2$.
\end{itemize}
\end{prop}
\proof
(i) By Proposition \ref{deformation0}
two admissible deformations must coincide for small $t$.
Then by the properties (iii-v) they have to be the same for
all $t\in I(t_0)$.

(ii)   Since
$f^{t_1,\tau_1}(b\D)\cap f^{t,\tau_2}(b\D)=\emptyset$
for $0\le t\le t_2$, then the intersection index of
$f^{t_1,\tau_1}(\bar\D)$ and $f^{t,\tau_2}(\bar\D)$
is independent of $t$ (see \cite{McS,MiWh}).
Since $f^{t_1,\tau_1}(\D)\cap f^{0,\tau_2}(\D)=\emptyset$,
then $f^{t_1,\tau_1}(\D)\cap f^{t_2,\tau_2}(\D)=\emptyset$,
also.
Q.E.D.
\medskip

An admissible deformation defined on a closed
interval $[0,t_0]$ can be extended to a larger interval.
We include a slightly stronger version of that result.

\begin{prop}
\label{IFT}
Let $(f^{t,\tau})$ be an admissible deformation defined
on $I(t_0) = [0,t_0)$. Suppose that for every
$\tau \in\R/2\pi \Z$ there exists an increasing sequence $(t^k)$,
$k = 1,2,...$ with $t^k \to t_0$ such that $f^{t^k,\tau}$
converges in the $C^{m}(\bar\D)$-norm for every $m$ to a
$J$-complex disc $f^{\infty,\tau}$. Then the deformation
can be extended to $I(t_1)$ for some $t_1>t_0$.
\end{prop}
\proof
It follows from Definition \ref{def-def}(iii)
that $f^{\infty,\tau}|_{b\D}$ is an embedding.
Since all  $f^{t,\tau}$ are embeddings, then
by the positivity of intersection indices
\cite{McS,MiWh}, the limiting disc
$f^{\infty,\tau}$ remains an embedding.
We obtain the discs $f^{t,\tau}$ for $t>t_0$ as
small deformations of the discs $f^{\infty,\tau}$
by applying the results \cite{GaSu,HoLiSi} about
small deformations of $J$-complex discs attached
to totally real manifolds. In the case of complex
dimension 2, such deformations are
governed by the normal Maslov index 
(see, e. g.,\cite{GaSu,HoLiSi}). The latter
is equal to zero for the disc $f^{\infty,\tau}$
because of the winding number condition in
Definition \ref{def-def}(iii).
Hence, the family extends to an interval $I(t_1)$
for some $t_1>t_0$.

It remains to show that the extended family on $I(t_1)$
still satisfies Definition \ref{def-def}, especially
part (iii).
For simplicity of notation we temporarily omit $\tau$
in $f^{t,\tau}$ and write just $f^t$.
The domain $\Omega$ naturally splits into two parts:
$\Omega = \Omega_1 \cup \Omega_2$, where
$\Omega_2 = H(\A_\delta \times \D)$ and
$\Omega_1 = \Omega \backslash \Omega_2$.
Likewise, every disc $D^t= f^t(\D)$
also splits into two parts:
$D^t= D^t_1 \cup D^t_2$.
Here $D^t_2$ is the subset of $\Omega_2$
represented as the graph of the function
$w^t:=w^{t,\tau}$ over ${\cal A}_\delta$,
and $D^t_1=D^t\backslash D^t_2$.
By Definition \ref{def-def}(iii), for every $0<t<t_0$,
we have $D^t_j = D^t\cap \Omega_j$, $j=1,2$.
We claim that the latter still holds for every $0<t<t_1$,
in particular, $D^t_1 \subset \Omega_1$.
Denote by $t'$ the supremum of the set of all $t<t_1$
for which $D^t_j = D^t\cap \Omega_j$, $j=1,2$.
Consider the Levi-flat hypersurface
$\Pi_\delta:= H(\{ |z|= \delta \}\times \D)$.
Note that by the Hopf lemma $D^t$ is transverse to
$\Pi_\delta$ for $t < t'$ (see, e. g., \cite{GaSu}).
If $t'< t_1$ then $D^{t'}_1$ touches $\Pi_\delta$ at
an interior point which is impossible
(see, e. g., \cite{DiSu}). Hence $t' = t_1$,
$D^t_1 \subset \Omega_1$ for all $0<t < t_1$,
and (iii) follows. The other conditions in
Definition \ref{def-def} are fulfilled automatically.
Q.E.D.

\section{Proof of the main result}

We establish {\it a priori} estimates for
any admissible deformation on an open interval.

\begin{lemma}
\label{Beltrami-mod1}
\begin{itemize}
\item[(i)]
Let $q$ and $Q$ be bounded functions in the annulus
$\A_\delta$ ($0<\delta<1$), $|q|\le q_0<1$, $|Q|\le Q_0$, here
$q_0$ and $Q_0$ are constants.
Let $\epsilon>0$ and $0<\delta'<\delta$.
Let $w$ be a solution of
\begin{equation}
\label{Beltrami-v}
w_{\bar z}=q w_z + Q.
\end{equation}
in $\A_\delta$,
such that $\epsilon\le |w| \le 1/\epsilon$ and $|w(z)|=1$
for $z\in b\D$.
Then $||w||_{C^\alpha(\A_{\delta'})}\le C$; here
$0<\alpha<1$ and $C>0$ depend on
$\epsilon$, $\delta$, $\delta'$, $q_0$, and $Q_0$
only.

\item[(ii)] Suppose in addition
that
$||q||_{C^{k,\beta}(\A_\delta)}
+||Q||_{C^{k,\beta}(\A_\delta)}\le Q_0$,
for some $0<\beta<1$ and $k\ge 0$. Then
$||w||_{C^{k+1,\beta}(\A_{\delta'})}\le C$; here
$C>0$ depends on
$\beta$, $k$, $\epsilon$, $\delta$, $\delta'$, $q_0$, and $Q_0$
only.
\end{itemize}
\end{lemma}

\proof
(i)
We apply the reflection principle.
Let $\A_\delta^*=\{ z^*: z\in\A_\delta \}$, here
$z^*:=1/\bar z$.
Extend $w$ and the coefficients $q$ and $Q$
to $\A_\delta^*$ by putting
$$
w(z)=(w(z^*))^*,\quad
q(z) = \bar {q(z^*)}z^2 (z^*)^2,\quad
Q(z) = \bar {Q(z^*)}(z^*)^2 w(z)^2
$$
for $z\in\A_\delta^*$. Then $w$ is continuous in
$G=\A_\delta\cup\A_\delta^*$ and satisfies
(\ref{Beltrami-v}) there.
The coefficients satisfy $|q|\le q_0$ and
$|Q|\le Q_0/\epsilon^2$ in $G$.

We claim that $w$ is uniformly bounded in $C^\alpha$
for some $0<\alpha<1$ after shrinking $G$.
It suffices to prove this fact for $G=\D$.
There exists a particular solution $w_0\in C^\alpha(\D)$
of the non-homogeneous equation (\ref{Beltrami-v}),
say by Proposition 2.1 (i) in \cite{SuTu2}.
Then $w=w_0+v$, where $v$ is a solution of the
homogeneous equation
$v_z=qv_{\bar z}$.
Then $v(z)=\phi(\xi(z))$, where
$\xi:\bar\D\to\bar\D$ is a fixed Beltrami homeomorphism
of class $C^\alpha(\D)$, and $\phi$ is holomorphic in $\D$.
Since $\phi$ is bounded, then the derivative $\phi'$
is also bounded in any smaller disc.
Then $v(z)=\phi(\xi(z))$, whence $w=w_0+v$ is
bounded in $C^\alpha$ in a smaller disc.

(ii) For simplicity, we first assume that $w$ has a
continuous logarithm in $\A_\delta$, which will be the case
in our applications. Then
$w=e^u$, $\Re u|_{b\D}=0$, and $u$ satisfies in $\A_\delta$
the equation
$$
u_{\bar z}=qu_z+Qw^{-1}.
$$
Let $\chi$ be a smooth cut-off function on $\bar\D$ vanishing in
a neighborhood of
$(1-\delta)\bar\D$ and such that $\chi\equiv 1$ in $\A_{\delta'}$.
Put $v=\chi u$. We can assume that $q$ and $Q$ are extended
over all of $\D$. Then
\begin{eqnarray}
\label{Beltrami-v1}
v_{\bar z}=qv_z+\chi Qw^{-1}+u(\chi_{\bar z}-q\chi_z)
\end{eqnarray}
By part (i), after shrinking $\delta$, we can assume
$w^{-1}, u \in C^\alpha(\A_{\delta})$.
Then $v$ satisfies (\ref{Beltrami-v1}) in $\D$
with boundary condition $\Re v|_{b\D}=0$.
Without loss of generality $v(1)=0$.
The conclusion now follows by successively
applying Proposition 2.1(ii) from \cite{SuTu2}.

In the general case we can use the above argument with
appropriate cut-off function to prove the result for every
proper sector of the annulus $\A_\delta$. Since the whole
annulus is a union of two such sectors, then the conclusion
will hold for $\A_\delta$.
Q.E.D.
\medskip

For an admissible deformation $(f^{t,\tau})$ we obtain
{\it a priori} estimates of the derivatives of the functions
$w^{t,\tau}$ from Definition \ref{def-def}(iii).
\begin{lemma}
\label{apriori}
Let $0<t_0<1$. Let $(f^{t,\tau})$ be an admissible deformation
on $0\le t<t_0$. Then for every $0<\delta'<\delta$ and $m\ge1$
there exists $C>0$ such that for every $0\le t<t_0$
and $\tau\in\R/2\pi\Z$ we have
$||w^{t,\tau}||_{C^m(\A(\delta'))}\le C$.
\end{lemma}
In fact the constant $C$ is independent of $t_0$ but we do not need
it in our application.
\medskip

\proof
Fix $t_1 < t_0$, say, $t_1 = t_0/2$.
For $0<t<t_1$ the desired estimate holds.
We need to show that it also holds for $t_1\le t<t_0$.
From the properties of the admissible deformation
it follows that the union of the discs $f^{t,\tau}$,
$0\le t<t_1$, $\tau\in\R/2\pi\Z$, cover an open
neighborhood of $f^0(\bar\D)$ in $H(\A_\delta\times\bar\D)$.
Since the discs $f^{t,\tau}$ do not intersect, then
for $t_1\le t<t_0$ the functions $w^{t,\tau}$ are uniformly
separated from zero.
Hence the desired conclusion follows by
Lemma \ref{Beltrami-mod1}. Q.E.D.
\medskip

The given symplectic form $\omega$ and almost complex
structure $J$ tamed by $\omega$ define a Riemannian metric
$$
\mu(V,W) = \frac{1}{2}(\omega(V,JW) + \omega(W,JV)).
$$
Let $f$ be a $J$-complex disc in $\bar\Omega$.
Let $E(f)$ denote the area of $f$ with respect
to $\mu$. Then (see, e. g., \cite{McS})
\begin{equation}
\label{def_area}
E(f) = \int_\D f^*\omega.
\end{equation}
We denote by $L(f)$ the length of the boundary of $f$, that is,
$$
L(f) = \int_0^{2\pi}
\left |\frac{df(e^{i\theta})}{d\theta}\right |_\mu d\theta,
$$
here $|\bullet|_\mu$ is the norm  defined by $\mu$.
Since the form $\omega$ is exact in a neighborhood of
$\bar\Omega$, that is, $\omega = d\lambda$, then by Stokes'
formula
\begin{equation}
\label{iso}
E(f) = \int_{f(\D)} \omega = \int_{f(b\D)} \lambda \leq C L(f),
\end{equation}
where $C > 0$ depends only on $\Omega$, $\lambda$, and $\mu$.
This inequality is a special case of the isoperimetric inequality
for $J$-complex curves, see \cite{Gr, McS}.

By Lemma \ref{apriori} the lengths of boundaries of $f^{t,\tau}$
are uniformly bounded. Hence we obtain an upper bound on areas
of the discs from an admissible deformation.
\begin{cor}
\label{Ascoli}
Let $0<t_0<1$. Let $(f^{t,\tau})$ be an admissible deformation
on $0\le t<t_0$. Then there exists a constant $C > 0$ such that
$E(f^{t,\tau}) \leq C$
for all $t$ and $\tau$.
\end{cor}
\medskip

\noindent
{\bf Proof of Theorem \ref{tori1}.}
Let $(f^{t,\tau})$ be an admissible deformation on
$[0, t_0)$. Consider a sequence $t^k \to t_0$ as
$k \to \infty$. Consider the sequence $(f^{t^k,\tau})$
for a fixed value of $\tau$. For simplicity of notation
we again omit $\tau$ in $f^{t,\tau}$ and write just $f^t$.

Since the areas of all discs $f^t$ are bounded, then
by Gromov's \cite{Gr} compactness theorem
(see also \cite{IvSu, McS}),
after passing to a subsequence if necessary, the sequence
$f^{t^k}$ converges to a $J$-complex disc $f^\infty$
uniformly on every compact subset of
$\overline\D \backslash \Sigma$. Here $\Sigma$ is a
finite set, where bubbles arise. The map $f^\infty$ is
smooth on $\overline\D$ and $f^\infty(b\D) \subset \Lambda^{t_0}$.
A bubble is a non-constant $J$-complex sphere
(a non-constant $J$-complex map from the Riemann sphere to $M$)
or a non-constant $J$-complex disc with boundary in $\Lambda^{t_0}$; disc-bubbles arise only
at the boundary points of $\D$.
We will prove that $\Sigma=\emptyset$, that is,
there are no bubbles. Then Gromov's compactness
theorem will imply the convergence in every $C^m(\bar\D)$-norm.

Since the form $\omega$ is exact, then by Stokes' formula every
$J$-complex sphere in $\Omega$ has zero area. Hence, there are
no spherical bubbles, and  $\Sigma$ can contain only points
of $b\D$, where disc-bubbles arise.
The sequence $f^{t^k}(\D)$ converges to a finite union
of $f^\infty(\D)$ and disc-bubbles in the Hausdorf
metric. It follows from the normalization condition of
Definition \ref{def-def}(iv) that the disc $f^\infty(\D)$
does not degenerate to a single point.

Let $F$ be the limit of the sequence $D^{t^k}$
in the Hausdorf metric.
Recall the decomposition $\Omega = \Omega_1 \cup \Omega_2$
that we used in the proof of Proposition \ref{IFT}.
Since $D^t_1 \subset \Omega_1$ for all $0<t < t_0$,
then $F \cap \Omega_2$ coincides with the Hausdorf limit of
the sequence $D^{t^k}_2$.
By Lemma \ref{apriori} and Ascoli's theorem
(after passing to a subsequence if necessary)
the sequence $w^{t^k}$
converges in $C^{m}(\A_\delta)$, $m \geq 0$, to some
function $w^\infty\in C^{m}(\A_\delta)$. But then
the graph of $w^\infty$ necessarily coincides with an open
piece of $f^\infty(\D)$. Hence boundary
bubbles do not arise either.

Thus, for every $\tau$ there is a subsequence
of $f^{t^k,\tau}$ converging in every $C^{m}$-norm.
Then by Proposition \ref{IFT} the admissible deformation
$(f^{t,\tau})$ extends past $t_0$. Hence, there is
an admissible deformation on the whole interval $[0,1]$,
and the proof of Theorem \ref{tori1} is complete.
Q.E.D.

\section{Non-squeezing}

We first establish the following
\begin{prop}
\label{tori2}
Let $G$ be a bounded domain in $\C^2$ and let $R > 0$.
Suppose $\bar G\subset R\D \times \C$ and
$\bar G\cap (R\D\times \{ 0 \})=\emptyset$.
Let $J$ be an almost complex structure on $\C^2$ tamed by
$\omega_\st$ and let $J=J_\st$ on $\C^2 \setminus\bar G$.
Then the domain  $R\D \times (\C\setminus\{ 0 \})$ is
foliated by a real one-parameter family of Levi-flat
hypersurfaces $\Gamma^t$, $t > 0$, with boundary
$b\Gamma^t=\Lambda^t = Rb\D \times tb\D$. Every
hypersurface $\Gamma^t$ in turn is foliated by embedded
$J$-complex discs attached to $\Lambda^t$, and
every such disc has Euclidean area $\pi R^2$.
\end{prop}

\proof
Without loss of generality assume $R = 1$.
Consider the Euclidean ball $s\B$ in $\C^2$ for $s > 0$
big enough so that $\bar G \subset s\B$. The boundary
of $s\B$ is a strictly pseudoconvex hypersurface with
respect to $J$. We apply Theorem \ref{tori1} to the family of tori
$\Lambda^t = b\D \times t b \D$ in $\Omega = s\B$.
As a result, we obtain a foliation by hypersurfaces $\Gamma^t$.
Every hypersurface $\Gamma^t$ in turn is foliated by embedded
$J$-complex discs $f^{t,\tau} = (z^{t,\tau},w^{t,\tau})$.

Since $J=J_\st$ on $\C^2 \backslash \bar G$, then
for big $t$, we have $f^{t,\tau}(\zeta)=(\zeta,t e^{i\tau})$.
Hence the hypersurfaces $\Gamma^t$ cover the set
$\D\times(\C\setminus r\D)$, where r is large.
By continuity, they cover the whole set
$\D\times(\C\setminus\{ 0 \})$ as stated.

We now claim that the discs $f^{t,\tau}$ have area $\pi$.
Indeed, let $\lambda = (i/2)(z\,d\bar z + w\,d\bar w)$.
Then $\omega_\st = d\lambda$.
Consider the parametrization of $\Lambda^t$ given by
$z = e^{i\phi}$, $w = t e^{i\psi}$. Then the restriction
$\theta = \lambda\vert_{\Lambda^t}$ has the form
$$
\theta = {1\over 2}\,(d\phi + t^2d\psi).
$$
By Stokes' formula (\ref{iso})
$$
E(f^{t,\tau}) = \int_{b\D} f^*\theta
 = {1\over 2} \int_{b\D} (z^{t,\tau})^*(d\phi)
 + {t^2\over 2} \int_{b\D} (w^{t,\tau})^*(d\psi) = \pi
$$
since the winding numbers of $z^{t,\tau}$ and $w^{t,\tau}$
are equal to $1$ and $0$ respectively for all $t$ and $\tau$.
Q.E.D.
\medskip

In particular, we immediately obtain the following
\begin{cor}
\label{tori3}
Let $J$ be a smooth almost complex structure in $R\D\times\C$
tamed by $\omega_\st$ and
such that $J-J_\st$ has compact support in $R\D\times\C$.
Then for every $p\in R\D\times\C$ there exists a $J$-complex
disc $f:\bar\D\to R\bar\D\times\C$ such that $f(0)=p$,
$f(b\D)\subset R\,b\D\times\C$, and $E(f)=\pi R^2$.
\end{cor}
\medskip

\noindent
{\bf Proof of Theorem \ref{rigidity}}.
Pushing forward the standard complex structure $J_\st$ by
$\phi$ yields an almost complex structure $J=\phi_*(J_\st)$
on $G_2$ tamed by $\omega_\st$.
Consider an exhaustion sequence of subdomains
$K_n\subset G_1$ such that every $K_n$ is relatively compact
in $K_{n+1}$.
Note that an almost complex structure $J$ is tamed
by $\omega_\st$ if and only if its complex matrix $A$
at every point has Euclidean norm $||A||<1$;
the set of such matrices is convex.
Therefore, for every $n$, there exists an almost complex
structure $\tilde J$ on $R\D \times \C$, which is tamed by
$\omega_\st$, coincides with $J$ on $\phi(K_n)$, and
coincides with $J_\st$ outside $\phi(K_{n+1})$.

Let $p = \phi(0)$. By Corollary \ref{tori3} there exists a
proper $\tilde J$-complex disc $D:=f(\D)$ in $R\D\times\C$
passing through $p$ with $E(D)=\pi R^2$.
Therefore $E(\phi(K_n)\cap D)\leq \pi R^2$.
Let $X_n= \phi^{-1}(D) \cap K_n$.
Since the map $\phi$ is a symplectomorphism,
then $E(X_n)\leq \pi R^2$.
Since $\phi: K_n \to \phi(K_n)$ is biholomorphic with respect to
$J_\st$ and $\tilde J$, then $X_n\in{\cal O}^1_0(K_n)$.
Since $(K_n)$ is an exhaustion sequence for $G_1$,
then by Bishop's convergence theorem (see, e. g., \cite{Ch}),
there is a subsequence of $(X_n)$ converging to a set
$X\in{\cal O}^1_0(G_1)$ with $E(X)\le\pi R^2$.
Hence, $\rh(G_1)\leq R$. Q.E.D.
\bigskip

\noindent
{\bf Proof of Corollary \ref{bidisc1}}.
By Theorem \ref{rigidity} it suffices to prove
that $\rh(\D^2_\R)>1$.
By Bishop's convergence theorem there exists
$X \in {\cal O}^1_0(\D^2_\R)$ such that
$E(X)=\pi(\rh(\D^2_\R))^2$.
The Euclidean unit ball $\B$ is contained in $\D^2_\R$ and
their boundaries $b\B$ and $b\D^2_\R$ meet
along two circles:
$$
S_1 = \{(z_1,z_2):x_1^2+x_2^2=1,\,y_1=y_2= 0 \},\quad
S_2 = \{(z_1,z_2):x_1 = x_2 = 0,\,y_1^2 + y_2^2 = 1 \}.
$$
Suppose that the boundary $bX = \overline X \backslash X$ of $X$
is contained in $S_1 \cup S_2$. Then $X$ is a complex
one-dimensional analytic subset in $\C^2\backslash(S_1\cup S_2)$.
By the reflection principle for analytic sets \cite{Ch},
$X$ extends as a complex $1$-dimensional analytic set to
a neighborhood of $S_1 \cup S_2$. Then by the uniqueness
theorem $X$ is contained in the complex algebraic curve
$(z^2_1+z^2_2)^2 =1$. But the latter does not pass through
the origin, a contradiction. Therefore, the closure $\overline X$
intersects the sphere $b\B$ at a point $p$ which is not in
$S_1 \cup S_2$. Since $X$ is closed in $\D^2_R$, the point $p$
is an interior point for $X$. The unit sphere $b\B$ is a strictly
pseudoconvex hypersurface. By the maximum principle \cite{Ch}
applied to the plurisubharmonic function $|z_1|^2+|z_2|^2$ on
the analytic set $X$, the latter can not be contained
in $\B$. Hence $X$ contains an open piece outside
the closed ball $\overline\B$.
Hence $E(X)>E(X \cap \B) \geq \pi$, and $\rh(\D^2_\R)>1$
as desired. Q.E.D.

\end{document}